\documentclass[11pt]{article}%
\usepackage{amsmath}
\usepackage{amsfonts}
\usepackage{amssymb}
\usepackage{graphicx}%
\newtheorem{theorem}{Theorem}

\newtheorem{definition}[theorem]{Definition}

\newtheorem{lemma}[theorem]{Lemma}

\newenvironment{proof}[1][Proof]{\noindent\textbf{#1.} }{\ \rule{0.5em}{0.5em}}
\begin{document}

\title{The smallest singular values of the icosahedral group}
\author{Charles F. Dunkl\thanks{Dept. of Mathematics, University of Virginia,
Charlottesville VA, 22904-4137; email: cfd5z@virginia.edu}}
\date{6 September 2018}
\maketitle

\begin{abstract}
For any finite reflection group $W$ on $\mathbb{R}^{N}$ and any irreducible
$W$-module $V$ there is a space of polynomials on $\mathbb{R}^{N}$ with values
in $V$. There are Dunkl operators parametrized by a multiplicity function,
that is, parameters asscociated with each conjugacy class of reflections. For
cerfain parameter values, called singular, there are nonconstant polynomials
annihilated by each Dunkl operator. There is a Gaussian bilinear form on the
polynomials which is positive for an open set of parameter values containing
the origin. When $W$ has just one class of reflections and $\dim V>1$  this
set is an interval bounded by the positive and negative singular values of
respective smallest absolute value. This interval is always symmetric around
$0$ for the symmetric groups. This property does not hold in general, and the
icosahedral group $H_{3}$ provides a counterexample. The interval for
positivity of the Gaussian form is determined for each of the ten irreducible
representations of $H_{3}$.

\end{abstract}

\section{Introduction}

Suppose $W$ is the finite reflection group generated by the reflections in the
reduced root system $R$ and $W$ has just one conjugacy class of reflections.
This means $R$ is a finite set of nonzero vectors in $\mathbb{R}^{N}$ such
that \thinspace$u,v\in R$ implies $\mathbb{R}u\cap R=\left\{  \pm u\right\}  $
and $v\sigma_{u}\in R$ where $\sigma_{u}$ is the reflection $x\mapsto
x-2\frac{\left\langle x,v\right\rangle }{\left\langle v,v\right\rangle }v$ and
$\left\langle \cdot,\cdot\right\rangle $ is the standard inner product.
Furthermore $R$ consists of one $W$-orbit. For a fixed vector $b_{0}$ such
that $\left\langle u,b_{0}\right\rangle \neq0$ for all $u\in R$ there is the
decomposition $R=R_{+}\cup R_{-}$ with $R_{+}:=\left\{  u\in R:\left\langle
u,b_{0}\right\rangle >0\right\}  $. The group $W$ is represented on the space
$\mathcal{P}$ of polynomials in $x=\left(  x_{1},\ldots,x_{N}\right)  $ by
$wp\left(  x\right)  =p\left(  xw\right)  $ for $w\in W$. Denote
$\mathbb{N}_{0}:=\left\{  0,1,2,\ldots\right\}  $ and for $\alpha\in
\mathbb{N}_{0}^{N}$ let $\left\vert \alpha\right\vert :=\sum_{i=1}^{N}%
\alpha_{i}$ and $x^{\alpha}:=\prod_{i=1}^{N}x_{i}^{\alpha_{i}}$, a monomial.
Then $\mathcal{P}=\mathrm{span}\left\{  x^{\alpha}:\alpha\in\mathbb{N}_{0}%
^{N}\right\}  $ and $\mathcal{P}_{n}=\mathrm{span}\left\{  x^{\alpha}%
:\alpha\in\mathbb{N}_{0}^{N},\left\vert \alpha\right\vert =n\right\}  $ the
space of polynomials homogeneous of degree $n$. Let $\kappa$ be a parameter
(called multiplicity function).

Suppose $\tau$ is an irreducible (unitary/orthogonal) representation of $W$ on
a (finite-dimensional) real vector space $V$ with basis $\left\{  u_{i}:1\leq
i\leq\dim V\right\}  $. The space $\mathcal{P}_{\tau}=\mathcal{P\otimes}V$ has
the basis $\left\{  x^{\alpha}\otimes u_{i}:\alpha\in\mathbb{N}_{0}^{N},1\leq
i\leq\dim V\right\}  $. There is a representation of $W$ on $\mathcal{P}%
_{\tau}$ defined to be the linear extension of%
\[
w\mapsto w\left(  p\left(  x\right)  \otimes u\right)  :=p\left(  xw\right)
\otimes\left(  \tau\left(  w\right)  u\right)  ,p\in\mathcal{P},u\in V.
\]
The associated Dunkl operators are the linear extension of%
\[
\mathcal{D}_{i}\left(  p\left(  x\right)  \otimes u\right)  =\frac{\partial
p\left(  x\right)  }{\partial x_{i}}\otimes u+\kappa\sum_{v\in R_{+}}%
\frac{p\left(  x\right)  -p\left(  x\sigma_{v}\right)  }{\left\langle
x,v\right\rangle }v_{i}\otimes\left(  \tau\left(  \sigma_{v}\right)  u\right)
.
\]
These operators mutually commute. If $\kappa=\kappa_{0}\in\mathbb{Q}$ and
there exists a nonconstant $p\in\mathcal{P}_{\tau}$ such that $\mathcal{D}%
_{i}p=0$ for all $i$ then $p$ is called a singular polynomial and $\kappa_{0}$
is a singular value. There is a Gaussian symmetric bilinear form $\left\langle
\cdot,\cdot\right\rangle _{\kappa}$ defined on $\mathcal{P}_{\tau}$, such that
$\left\langle wf,wg\right\rangle _{\kappa}=\left\langle f,g\right\rangle
_{\kappa}$ and $\left\langle \mathcal{D}_{i}f,g\right\rangle _{\kappa
}=\left\langle f,\left(  x_{i}-\mathcal{D}_{i}\right)  g\right\rangle
_{\kappa}$ for all $f,g\in\mathcal{P}_{\tau}$, $w\in W$ and $1\leq i\leq N$.
Shelley-Abrahamson \cite{SA} proved there is an interval for $\kappa$
containing zero for which the bilinear form is realized as an integral with
respect to a positive matrix-valued measure on $\mathbb{R}^{N}$. In \cite{D1}
the author showed that this interval is bounded by the smallest (in absolute
value) singular values. When $W$ is a symmetric group (type $A_{N-1}$) the
interval is always symmetric about $\kappa=0$. It was conjectured that this
always holds when $\dim V_{\tau}\geq2$. It is the purpose of this note to show
that the group $H_{3}$ provides a counterexample to the conjecture. We will
determine the smallest singular values for each of the ten (conjugacy classes
of) irreducible representations of $H_{3}$. 

The paper \cite{DJO} showed that the singular values for Coxeter groups with
one class of reflections are of the form $m/d_{j}$ where the fundamental
degrees of $W$ are $\left\{  d_{1},d_{2},\ldots,d_{N}\right\}  $,
$m\in\mathbb{Z}$ and $m/d_{j}\notin\mathbb{Z}$. Etingof and Stoica \cite{ES}
analyzed unitary representations of the rational Cherednik algebra of $W$ (the
abstract algebra generated by $\left\{  \mathcal{D}_{i}\right\}  $,
multiplication by $x_{i}$ ($1\leq i\leq N$) and translation by $w\in W$); the
parameter intervals allowing the unitary property include the intervals for
positivity of the Gaussian inner product. Balagovic and Puranik \cite{BP}
studied the representations of the Cherednik algebra of $H_{3}$; this involved
a much more detailed analysis than the singular values appearing here.

\subsection{General facts:}

Suppose that $\rho$ is an irreducible representation of $W$ with character
$\chi_{\rho}$ and suppose that $p\in\mathcal{P}_{\tau}$ has the property that
$\mathrm{span}\left\{  wp:w\in W\right\}  $ is a $W$-submodule on which $W$
acts by $\rho$, then $p$ is said to be of \textit{isotype} $\rho$. This
implies $\sum\limits_{v\in R_{+}}\sigma_{v}p=\varepsilon\left(  \rho\right)
p$ where $\varepsilon\left(  \rho\right)  =\#R_{+}\dfrac{\chi_{\rho}\left(
\sigma_{v}\right)  }{\chi_{\rho}\left(  1\right)  }$, the scalar by which the
central element $\sum\limits_{v\in R_{+}}\sigma_{v}$ of $\mathbb{R}W$ acts on
the module associated with $\rho$. These values are of key importance in the
analysis of singular polynomials. The operators $\mathcal{D}_{i}$ transform
under $W$ as follows: for $a\in\mathbb{R}^{N}$ define $\mathcal{D}_{a}%
p:=\sum_{i=1}^{N}a_{i}\mathcal{D}_{i}p$, then it can be shown that
$\mathcal{D}_{a}wp=w\mathcal{D}_{aw}p$ (this relies on the fact that $w$
permutes the elements of $R$). As a consequence any singular polynomial can be
expanded as a sum of singular polynomials with each of some isotype (the
decomposability of a representation of $W$).

\begin{lemma}
\label{sgval}Suppose $p\in\mathcal{P}_{\tau}$ is homogeneous of degree $n$, of
isotype $\rho$, and singular for $\kappa=\kappa_{0}$ then $\kappa
_{0}=n/\left(  \varepsilon\left(  \rho\right)  -\varepsilon\left(
\tau\right)  \right)  $.
\end{lemma}

\begin{proof}
Apply the operator $\sum_{i=1}^{N}x_{i}\mathcal{D}_{i}$ to a generic
polynomial in $\mathcal{P}_{\tau}$ of the form $p=\sum_{j=1}^{m}p_{j}\left(
x\right)  \otimes u_{j}$ to obtain%
\begin{align*}
\sum_{i=1}^{N}x_{i}\mathcal{D}_{i}p  & =\sum_{i=1}^{N}x_{i}\frac{\partial
}{\partial x_{i}}p+\kappa\sum_{j=1}^{m}\sum_{v\in R_{+}}\left(  p_{j}\left(
x\right)  \otimes\tau\left(  \sigma_{v}\right)  u_{j}-p_{j}\left(  x\sigma
_{v}\right)  \otimes\tau\left(  \sigma_{v}\right)  u_{j}\right)  \\
& =np+\kappa\left\{  \varepsilon\left(  \tau\right)  \sum_{j=1}^{m}%
p_{j}\left(  x\right)  \otimes u_{j}-\sum_{v\in R_{=}}\sigma_{v}\left(
p_{j}\left(  x\right)  \otimes u_{j}\right)  \right\}  \\
& =np+\kappa\left(  \varepsilon\left(  \tau\right)  -\varepsilon\left(
\rho\right)  \right)  p.
\end{align*}
Thus $\sum_{i=1}^{N}x_{i}\mathcal{D}_{i}p=0$ implies $n+\kappa\left(
\varepsilon\left(  \tau\right)  -\varepsilon\left(  \rho\right)  \right)  =0$.
\end{proof}

Let $\det$ denote the sign-character of $W$, then $\tau^{\vee}:w\mapsto
\det\left(  w\right)  \tau\left(  w\right)  $ is an irreducible representation
and the character $\chi_{\tau^{\vee}}\left(  w\right)  =\det\left(  w\right)
\chi_{\tau}\left(  w\right)  $. This induces an identification between
$V_{\tau}$ and $V_{\tau^{\vee}}$ so that $\tau^{\vee}\left(  w\right)
u=\det\left(  w\right)  \tau\left(  w\right)  u$ for $u\in V_{\tau}$. Lift
this map to $\mathcal{P}_{\tau}\rightarrow\mathcal{P}_{\tau^{\vee}}$. Since
$\det\left(  \sigma_{v}\right)  =-1$ we see that a singular polynomial $p$ in
$\mathcal{P}_{\tau}$ for the singular value $\kappa_{0}$ is mapped to a
singular polynomial $p^{\vee}$ in $\mathcal{P}_{\tau^{\vee}}$ for $-\kappa
_{0}$; and if $p$ is of isotype $\rho$ then $p^{\vee}$ is of isotype
$\rho^{\vee}$. We will refer to this as the \textit{duality} principle.

\section{The icosahedral group}

This is the symmetry group of the regular icosahedron%
\[
\mathcal{Q}_{12}=\left\{  \left(  0,\pm\tau,\pm1\right)  ,\left(  \pm
1,0,\pm\tau\right)  ,\left(  \pm\tau,\pm1,0\right)  \right\}
\]
(12 vertices, 20 triangular faces with 5 meeting at each vertex, 30 edges) and
of the dual dodecahedron $\mathcal{Q}_{20}$ (vertices at the mid-points of the
faces of $Q_{20}$, 12 pentagonal faces with 3 meeting at each vertex), where
$\tau:=\frac{1}{2}\left(  1+\sqrt{5}\right)  $ and $\tau^{2}=\tau+1$. The root
system $R$ of $H_{3}$  consists of the cyclic permutations of $\left(
\pm2,0,0\right)  $ and $\left(  \pm\tau,\pm\left(  \tau-1\right)
,\pm1\right)  $ (independent choices of signs). The positive roots $R_{+}$ are
those having positive inner product with $\left(  3,2\tau,1\right)  $ and
$\#R_{+}=15$.

The character table of $H_{3}$ is
\[%
\begin{vmatrix}
1\\
1\\
15\\
15\\
20\\
12\\
12\\
20\\
12\\
12
\end{vmatrix}
\circ%
\begin{vmatrix}
1 & 1 & 3 & 3 & 3 & 3 & 4 & 4 & 5 & 5\\
1 & -1 & -3 & 3 & -3 & 3 & 4 & -4 & 5 & -5\\
1 & -1 & 1 & -1 & 1 & -1 & 0 & 0 & 1 & -1\\
1 & 1 & -1 & -1 & -1 & -1 & 0 & 0 & 1 & 1\\
1 & 1 & 0 & 0 & 0 & 0 & 1 & 1 & -1 & -1\\
1 & 1 & 1-\tau & 1-\tau & \tau & \tau & -1 & -1 & 0 & 0\\
1 & 1 & \tau & \tau & 1-\tau & 1-\tau & -1 & -1 & 0 & 0\\
1 & -1 & 0 & 0 & 0 & 0 & 1 & -1 & -1 & 1\\
1 & -1 & -\tau & \tau & \tau-1 & 1-\tau & -1 & 1 & 0 & 0\\
1 & -1 & \tau-1 & 1-\tau & -\tau & \tau & -1 & 1 & 0 & 0
\end{vmatrix}
\]
The column on the left lists the cardinalities $n_{i}$ of the conjugacy
classes $C_{i}$. In this list $C_{1}=\left\{  I\right\}  ,C_{2}=\left\{
-I\right\}  $ and $C_{3}$ is the class of reflections. There are three classes
of rotations: $C_{4}$ contains the period two rotations (with axes being the
mid-points of the edges of $Q_{12}$, $C_{5}$ contains the period three
rotations with axes at the dodecahedron vertices, and $C_{6},C_{7}$ contains
the period five rotations with axes at the icosahedron vertices, with angles
of $\frac{4\pi}{5},\frac{2\pi}{5}$ respectively. The other three classes come
from the product of $-I$ with $C_{4},C_{5},C_{6}$. For the purpose of
computing Poincar\'{e} series we list the polynomials $d_{i}\left(  z\right)
=\det\left(  I-zw\right)  $ where $w\in C_{i}$%
\[%
\begin{vmatrix}
C_{1} & \left(  1-z\right)  ^{3}\\
C_{2} & \left(  1+z\right)  ^{3}\\
C_{3} & \left(  1+z\right)  \left(  1-z\right)  ^{2}\\
C_{4} & \left(  1-z\right)  \left(  1+z\right)  ^{2}\\
C_{5} & \left(  1-z\right)  \left(  1+z+z^{2}\right)  \\
C_{6} & \left(  1-z\right)  \left(  1+\tau z+z^{2}\right)  \\
C_{7} & \left(  1-z\right)  \left(  1+\left(  1-\tau\right)  z+z^{2}\right)
\\
C_{8} & \left(  1+z\right)  \left(  1-z+z^{2}\right)  \\
C_{9} & \left(  1+z\right)  \left(  1+\left(  \tau-1\right)  z+z^{2}\right)
\\
C_{10} & \left(  1+z\right)  \left(  1-\tau z+z^{2}\right)  .
\end{vmatrix}
\]
Denote the character defined in column $\#i$ of the table by $\chi_{i}$ and
the corresponding representation by $\tau_{i}$ acting on the $\chi_{i}\left(
C_{1}\right)  $-dimensional space $V_{i}$. Thus $\tau_{1}=\mathrm{triv}$,
$\tau_{2}=\det$ and $\tau_{3}$ is the reflection representation. By the
Poincar\'{e} series argument the character for $H_{3}$ acting on
$\mathcal{P}_{n}$ evaluated at $C_{i}$ is the coefficient of $z^{n}$ in
$1/d_{i}\left(  z\right)  $. Thus the character of $H_{3}$ for $\mathcal{P}%
_{n}\otimes V_{i}$ at $C_{i}$ is the coefficient of $z^{n}$ in $\chi
_{i}\left(  C_{i}\right)  /d_{i}\left(  z\right)  $. The inner product of
characters is used to determine multiplicities.

\begin{definition}
For $1\leq i,j\leq10$ let%
\[
\boldsymbol{P}_{ij}\left(  z\right)  =\frac{1}{120}\sum_{s=1}^{10}n_{s}%
\frac{\chi_{i}\left(  C_{s}\right)  \chi_{j}\left(  C_{s}\right)  }%
{d_{s}\left(  z\right)  }%
\]
and let $L_{ij}$ denote the lowest power of $z$ appearing (with a nonzero
coefficient) in the power series expansion of $\boldsymbol{P}_{ij}\left(
z\right)  $.
\end{definition}

Of course $\boldsymbol{P}_{11}\left(  z\right)  =\dfrac{1}{\left(
1-z^{2}\right)  \left(  1-z^{6}\right)  \left(  1-z^{10}\right)  }$, the
Poincar\'{e} series for the $H_{3}$-invariant polynomials, and $2,6,10$ are
the fundamental degrees. Suppose $U$ is an $H_{3}$-invariant subspace of some
space $\mathcal{P}_{n}\otimes V_{i}$ on which $H_{3}$ acts isomorphically to
$\tau_{j}$ then $U$ is said to be of isotype $\tau_{j}$. Furthermore
$\mathcal{P}_{n}\otimes V_{i}$ has a complete decomposition into subspaces of
various isotypes. Thus the coefficient of $z^{n}$ in $\boldsymbol{P}%
_{ij}\left(  z\right)  $ is the number of components of $\mathcal{P}%
_{n}\otimes V_{i}$ which are of isotype $\tau_{j}$ (or components of
$\mathcal{P}_{n}\otimes V_{j}$ of isotype $\tau_{i}$).

The eigenvalues of $\sum\limits_{v\in R_{=}}\sigma_{v}$ acting on $V_{i}$ are
denoted $\varepsilon\left(  \tau_{i}\right)  =15\chi_{i}\left(  C_{3}\right)
/\chi_{i}\left(  C_{1}\right)  $ and are listed in $\left[
15,-15,5,-5,5,-5,0,0,3,-3\right]  $. From Lemma \ref{sgval} if there are
singular polynomials of isotype $\tau_{j}$ in $\mathcal{P}_{n}\otimes V_{i}$
then the singular value $\kappa_{0}=\dfrac{n}{\varepsilon\left(  \tau
_{j}\right)  -\varepsilon\left(  \tau_{i}\right)  }$. For each $i$ the list
$\left[  \dfrac{L_{ij}}{\varepsilon\left(  \tau_{j}\right)  -\varepsilon
\left(  \tau_{i}\right)  }:1\leq j\leq10,\varepsilon\left(  \tau_{j}\right)
\neq\varepsilon\left(  \tau_{i}\right)  \right]  $ contains the possible
singular values of minimal absolute values. That is, no singular value can
have a smaller absolute value than the extreme values of the list. There are
two parts of the argument: firstly to compute the lists for each $\tau_{i}$
(actually just $\tau_{1},\tau_{3},\tau_{5},\tau_{8},\tau_{9}$ because the
other cases follow from the duality principle), secondly to show that the
extreme values do correspond to singular polynomials. In each of the following
there is a list with 10 pairs $\left[  \left[  L_{ij},\varepsilon\left(
\tau_{j}\right)  -\varepsilon\left(  \tau_{i}\right)  \right]  \right]
_{j=1}^{10}$ and $\kappa_{-},\kappa_{+}$ denote the negative and positive
ratios of minimum absolute values; and each is given with the corresponding isotype

\begin{itemize}
\item $\tau_{1}$, $[0,0],[15,-30],[1,-10],[6,-20],[3,-10],[8,-20],[4,-15]$,
$[3,-15],$\newline$[2,-12],[5,-18]$; there is no $\kappa_{+}$ and $\kappa
_{-}=-\frac{1}{10}$ for $\tau_{3}$. The dual is $\tau_{2}$ with $\kappa
_{+}=\frac{1}{10}$ for $\tau_{4}$.

\item $\tau_{3}$, $[[1,10],[6,-20],[0,0],[1,-10],[2,0],[3,-10],[3,-5],[2,-5],$%
\newline$[1,-2],[2,-8]]$; $\kappa_{+}=\frac{1}{10}$ for $\tau_{1}$ and
$\kappa_{-}=-\frac{1}{10}$ for $\tau_{4}.$The dual is $\tau_{4}$ with
$\kappa_{-}=-\frac{1}{10}$ for $\tau_{2}$ and $\kappa_{+}=\frac{1}{10}$ for
$\tau_{3}$.

\item $\tau_{5}$, $[[3,10],[8,-20],[2,0],[3,-10],[0,0],[3,-10],[1,-5],[2,-5],$%
\newline$[1,-2],[2,-8]]$; $\kappa_{+}=\frac{3}{10}$ for $\tau_{1}$ and
$\kappa_{-}=-\frac{1}{5}$ for $\tau_{7}$. The dual is $\tau_{6}$ with
$\kappa_{+}=\frac{1}{5}$ for $\tau_{8}$ and $\kappa_{-}=-\frac{3}{10}$ for
$\tau_{2}$.

\item $\tau_{8}%
,~[[3,15],[4,-15],[2,5],[3,-5],[2,5],[1,-5],[1,0],[0,0],[1,3],[2,-3]]$;
$\kappa_{+}=\frac{1}{5}$ for $\tau_{1}$ and $\kappa_{-}=-\frac{1}{5}$ for
$\tau_{6}$. The dual is $\tau_{7}$ with $\kappa_{+}=\frac{1}{5}$ for $\tau
_{5}$ and $\kappa_{-}=-\frac{1}{5}$ for $\tau_{5}$.

\item $\tau_{9}$%
,$~[[2,12],[5,-18],[1,2],[2,-8],[1,2],[2,-8],[2,-3],[1,-3],[0,0],[1,-6]]$;
$\kappa_{+}=\frac{1}{6}$ for $\tau_{1}$ and $\kappa_{-}=-\frac{1}{6}$ for
$\tau_{10}$. The dual is $\tau_{10}$ with $\kappa_{+}=\frac{1}{6}$ for
$\tau_{9}$ and $\kappa_{-}=-\frac{1}{6}$ for $\tau_{2}$.
\end{itemize}

The key point here is the asymmetry in $\left(  \kappa_{-},\kappa_{+}\right)
$ manifested in $\tau_{5}$ and $\tau_{6}$.

\subsection{Construction of singular polynomials}

The idea is to set up a framework that can be run by symbolic computation
software. To take the place of the abstract $H_{3}$-modules $V_{i}$ we will
use submodules of polynomials in $t=\left(  t_{1},t_{2},t_{3}\right)  $. The
generic element of $\mathcal{P}_{\tau_{i}}$ is a polynomial $p\left(
x,t\right)  $ homogeneous in each of $x$ and $t$. The action of $H_{3}$ is
given by $wp\left(  x,t\right)  =p(xw,tw)$ and the Dunkl operators are%
\[
\mathcal{D}_{i}p\left(  x,t\right)  =\frac{\partial p\left(  x,t\right)
}{\partial x_{i}}+\kappa\sum_{v\in R_{+}}\frac{p\left(  x,t\sigma_{v}\right)
-p\left(  x\sigma_{v},t\sigma_{v}\right)  }{\left\langle x,v\right\rangle
}v_{i}.
\]
Actually it is more efficient to compute the gradient%
\[
\nabla\left(  \kappa\right)  p\left(  x,t\right)  :=\left[  \dfrac{\partial
p\left(  x,t\right)  }{\partial x_{i}}\right]  _{i=1}^{3}+\kappa
\sum\limits_{v\in R_{+}}\dfrac{p\left(  x,t\sigma_{v}\right)  -p\left(
x\sigma_{v},t\sigma_{v}\right)  }{\left\langle x,v\right\rangle }v
\]
so as to avoid repeated evaluations of the 15 difference quotients. To
identify the various components use the operators
\begin{align*}
Sp\left(  x,t\right)    & :=\sum_{v\in R_{+}}p\left(  x\sigma_{v},t\sigma
_{v}\right)  ,\\
S^{\left(  t\right)  }p\left(  x,t\right)    & :=\sum_{v\in R_{+}}p\left(
x,t\sigma_{v}\right)
\end{align*}
By use of the series $\boldsymbol{P}_{ij}\left(  z\right)  $ we can determine
the isotypes appearing in the lowest degree polynomials:%
\[
\mathcal{P}_{1}\cong V_{3},~\mathcal{P}_{2}\cong V_{1}\oplus V_{9}%
,~\mathcal{P}_{3}\cong V_{3}\oplus V_{5}\oplus V_{8}.
\]
The unwanted components $V_{1}$ of $\mathcal{P}_{2}$, and $V_{3}$ of
$\mathcal{P}_{3}$ (the first is $\sum_{i=1}^{3}x_{i}^{2}$ and the second is
$\sum_{i=1}^{3}x_{i}^{2}\mathcal{P}_{1}$) can be eliminated by restricting to
harmonic polynomials, that is, those in the kernel of $\sum_{i=1}%
^{3}\mathcal{D}_{i}^{2}$. Converting to the variable $t$ we obtain bases for
the component $V_{9}$ in $\mathcal{P}_{2}$ and for $V_{5}\oplus V_{8}$ in
$\mathcal{P}_{3}$
\begin{align*}
&  t_{1}^{2}-t_{3}^{2},t_{2}^{2}-t_{3}^{2},t_{1}t_{2},t_{1}t_{3},t_{2}t_{3};\\
&  t_{1}\left(  t_{2}^{2}-t_{3}^{2}\right)  ,t_{2}\left(  t_{3}^{2}-t_{1}%
^{2}\right)  ,t_{3}\left(  t_{1}^{2}-t_{2}^{2}\right)  ,t_{1}\left(  t_{1}%
^{2}-3t_{2}^{2}\right)  ,t_{2}\left(  t_{2}^{2}-3t_{3}^{2}\right)
,t_{3}\left(  t_{3}^{2}-3t_{1}^{2}\right)  ,t_{1}t_{2}t_{3}.
\end{align*}
By applying $S^{\left(  t\right)  }$ to the degree 3 polynomials and using the
eigenvalues $-5,0$ for $V_{5},V_{8}$ respectively we find a basis for the
$V_{5}$ and $V_{8}$ components:
\begin{align*}
V_{5}:  & t_{1}\left(  t_{1}^{2}-3\left(  1-\tau\right)  t_{2}^{2}-3\tau
t_{3}^{2}\right)  ,~t_{2}\left(  t_{2}^{2}-3\left(  1-\tau\right)  t_{3}%
^{2}-3\tau t_{1}^{2}\right)  ,\\
& t_{3}\left(  t_{3}^{2}-3\left(  1-\tau\right)  t_{1}^{2}-3\tau t_{2}%
^{2}\right)  ,
\end{align*}%
\begin{align*}
V_{8}  & :t_{1}\left(  t_{1}^{2}-\left(  1+\tau\right)  t_{2}^{2}+\left(
\tau-2\right)  t_{3}^{2}\right)  ,~t_{2}\left(  t_{2}^{2}-\left(
1+\tau\right)  t_{3}^{2}+\left(  \tau-2\right)  t_{1}^{2}\right)  ,~\\
& t_{3}\left(  t_{3}^{2}-\left(  1+\tau\right)  t_{1}^{2}+\left(
\tau-2\right)  t_{2}^{2}\right)  ,~t_{1}t_{2}t_{3}.
\end{align*}
Denote the entries in the list for $V_{5}$ by $p_{i}^{\left(  5\right)
}\left(  t\right)  ,i=1\ldots3$ respectively, and those in the list for
$V_{8}$ by $p_{i}^{\left(  8\right)  }\left(  t\right)  ,i=1\ldots4$
respectively. Our notation for singular polynomials is $h_{i}^{\left(
j\right)  }\left(  x,t\right)  $ which means that $h_{i}^{\left(  j\right)
}\left(  x,t\right)  \in\mathcal{P}\otimes V_{j}$ and $h_{i}^{\left(
j\right)  }\left(  x,t\right)  $ is of isotype $\tau_{i}$. Define%
\begin{align*}
h_{1}^{\left(  3\right)  }\left(  x,t\right)    & :=\sum_{i=1}^{3}x_{i}%
t_{i},\\
h_{1}^{\left(  9\right)  }\left(  x,t\right)    & :=\sum_{1\leq i<j\leq
3}\left\{  \left(  x_{i}^{2}-x_{j}^{2}\right)  \left(  t_{i}^{2}-t_{j}%
^{2}\right)  +6x_{i}x_{j}t_{i}t_{j}\right\}  ,\\
h_{1}^{\left(  5\right)  }\left(  x,t\right)    & :=\sum_{i=1}^{3}%
p_{i}^{\left(  5\right)  }\left(  x\right)  p_{i}^{\left(  5\right)  }\left(
t\right)  ,\\
h_{1}^{\left(  8\right)  }\left(  x,t\right)    & :=\sum_{i=1}^{3}%
p_{i}^{\left(  8\right)  }\left(  x\right)  p_{i}^{\left(  8\right)  }\left(
t\right)  +20p_{4}^{\left(  8\right)  }\left(  x\right)  p_{4}^{\left(
8\right)  }\left(  t\right)  .
\end{align*}
Then $Sh_{i}^{\left(  j\right)  }\left(  x,t\right)  =15h_{1}^{\left(
j\right)  }\left(  x,t\right)  $ for $j=3,5,8,9$ (that is, $h_{1}^{\left(
j\right)  }$ is of isotype $\tau_{1}=\mathrm{triv}$). By direct computation%
\begin{align*}
\nabla\left(  \frac{1}{10}\right)  h_{1}^{\left(  3\right)  }\left(
x,t\right)    & =0,~\nabla\left(  \frac{3}{10}\right)  h_{1}^{\left(
5\right)  }\left(  x,t\right)  =0,\\
\nabla\left(  \frac{1}{5}\right)  h_{1}^{\left(  8\right)  }\left(
x,t\right)    & =0,~\nabla\left(  \frac{1}{6}\right)  h_{1}^{\left(  9\right)
}\left(  x,t\right)  =0.
\end{align*}
It remains to construct the $\left(  \tau_{1},\kappa=-\frac{1}{10}\right)  $,
$\left(  \tau_{3},\kappa=-\frac{1}{10}\right)  $, $\left(  \tau_{5}%
,\kappa=-\frac{1}{5}\right)  $, $\left(  \tau_{8},\kappa=-\frac{1}{5}\right)
$, and $\left(  \tau_{9},\kappa=-\frac{1}{6}\right)  $ singular polynomial.
For $\tau_{1}$ the solution is simply $h_{3}^{\left(  1\right)  }\left(
x\right)  :=\sum_{i=1}^{3}c_{i}x_{i}$ with arbitrary constants, of isotype
$\tau_{3}$. The solution for $\tau_{3}$ is of first degree in $x$ and $t$:%
\begin{gather*}
h_{4}^{\left(  3\right)  }\left(  x,t\right)  :=c_{1}\left(  x_{2}t_{3}%
-x_{3}t_{2}\right)  +c_{2}\left(  x_{3}t_{1}-x_{1}t_{3}\right)  +c_{3}\left(
x_{1}t_{2}-x_{2}t_{1}\right)  ,\\
\nabla\left(  -\frac{1}{10}\right)  h_{4}^{\left(  3\right)  }\left(
x,t\right)  =0,~Sh_{4}^{\left(  3\right)  }\left(  x,t\right)  =-5h_{4}%
^{\left(  3\right)  }\left(  x,t\right)  .
\end{gather*}
For $\tau_{5}$ we find%
\begin{align*}
h_{7}^{\left(  5\right)  }\left(  x,t\right)    & :=c_{1}\left(  x_{1}%
p_{2}^{\left(  5\right)  }\left(  t\right)  +\left(  \tau+1\right)  x_{2}%
p_{3}^{\left(  5\right)  }\left(  t\right)  \right)  +c_{2}\left(  \left(
1+\tau\right)  x_{1}p_{1}^{\left(  5\right)  }\left(  t\right)  +x_{3}%
p_{3}^{\left(  5\right)  }\left(  t\right)  \right)  \\
& +c_{3}\left(  x_{2}p_{1}^{\left(  5\right)  }\left(  t\right)  +\left(
\tau+1\right)  x_{3}p_{2}^{\left(  5\right)  }\left(  t\right)  \right)
+c_{4}\left(  x_{3}p_{1}^{\left(  5\right)  }\left(  t\right)  +x_{2}%
p_{2}^{\left(  5\right)  }\left(  t\right)  +x_{1}p_{3}^{\left(  5\right)
}\left(  t\right)  \right)  ,
\end{align*}%
\[
\nabla\left(  -\frac{1}{5}\right)  h_{7}^{\left(  5\right)  }\left(
x,t\right)  =0,~Sh_{7}^{\left(  5\right)  }\left(  x,t\right)  =0.
\]
For $\tau_{8}$ the singular polynomials are%
\begin{align*}
h_{6}^{\left(  8\right)  }\left(  x,t\right)    & :=c_{1}\left(  \left(
3-\tau\right)  x_{1}p_{1}^{\left(  8\right)  }\left(  t\right)  +\left(
2+\tau\right)  x_{3}p_{3}^{\left(  8\right)  }\left(  t\right)  +10x_{2}%
p_{4}^{\left(  8\right)  }\left(  t\right)  \right)  \\
& +c_{2}\left(  \left(  2+\tau\right)  x_{2}p_{1}^{\left(  8\right)  }\left(
t\right)  +\left(  3-\tau\right)  x_{3}p_{2}^{\left(  8\right)  }\left(
t\right)  +10x_{1}p_{4}^{\left(  8\right)  }\left(  t\right)  \right)  \\
& +c_{3}\left(  \left(  2+\tau\right)  x_{1}p_{2}^{\left(  8\right)  }\left(
t\right)  +\left(  3-\tau\right)  x_{2}p_{3}^{\left(  8\right)  }\left(
t\right)  +10x_{3}p_{4}^{\left(  8\right)  }\left(  t\right)  \right)  ,
\end{align*}%
\[
\nabla\left(  -\frac{1}{5}\right)  h_{6}^{\left(  8\right)  }\left(
x,t\right)  =0,~Sh_{6}^{\left(  8\right)  }\left(  x,t\right)  =-5h_{6}%
^{\left(  8\right)  }\left(  x,t\right)  .
\]
Finally for $\tau_{9}$ we have%
\begin{align*}
h_{10}^{\left(  9\right)  }\left(  x,t\right)    & =c_{1}\left(  x_{3}\left(
t_{1}^{2}-t_{2}^{2}\right)  -\left(  x_{1}t_{1}-x_{2}t_{2}\right)
t_{3}\right)  +c_{2}\left(  x_{2}\left(  t_{1}^{2}-t_{3}^{2}\right)  -\left(
x_{1}t_{1}-x_{3}t_{3}\right)  t_{2}\right)  \\
& +c_{3}\left(  x_{1}\left(  t_{3}^{2}-t_{2}^{2}\right)  -\left(  x_{3}%
t_{3}-x_{2}t_{2}\right)  t_{1}\right)  +c_{4}\left(  x_{1}t_{2}-x_{2}%
t_{1}\right)  t_{3}+c_{5}\left(  x_{1}t_{3}-x_{3}t_{1}\right)  t_{2}.
\end{align*}%
\[
\nabla\left(  -\frac{1}{6}\right)  h_{10}^{\left(  9\right)  }\left(
x,t\right)  =0,~Sh_{10}^{\left(  9\right)  }\left(  x,t\right)  =-3h_{10}%
^{\left(  9\right)  }\left(  x,t\right)  .
\]
We have shown that the smallest singular values for $\tau_{5}$ and $\tau_{6}$
are $-\frac{1}{5},\frac{3}{10}$ and $-\frac{3}{10},\frac{1}{5}$ respectively.
For most cases the interval is symmetric but this is not a general principle.


\begin{thebibliography}{9}                                                                                                %
\bibitem {BP}M. Balagovic and A. Puranik, Irreducible representations of the
rational Cherednik algebra associated to the Coxeter group $H_{3}$, \textit{J.
Algebra }\textbf{405} (2014), 259-290, arxiv:1004.2108v3.

\bibitem {D1}C. Dunkl, The smallest singular values and vector-valued Jack
polynomials, arxiv:1804.09158v3.

\bibitem {DJO}C. Dunkl, M. de Jeu, and E. Opdam, Singular polynomials for
finite reflection groups, \textit{Trans. Amer. Math. Soc. }\textbf{346
}(1994), 237-256.

\bibitem {ES}P. Etingof and E. Stoica, Unitary representations of rational
Cherednik algebras, (with an appendix by S. Griffeth), \textit{Represent.
Theory} \textbf{13} (2009), 349-370, arXiv:0901.4595.

\bibitem {SA}S. Shelley-Abrahamson, The Dunkl weight function for rational
Cherednik algebras, 2018, arXiv:1803.00440.
\end{thebibliography}
\end{document}